\input amstex
\input amsppt.sty
\magnification=\magstep1
\hsize=30truecc
\vsize=22.2truecm
\baselineskip=16truept
\NoBlackBoxes
\TagsOnRight \pageno=1 \nologo
\def\Z{\Bbb Z}
\def\N{\Bbb N}

\def\l{\left}
\def\r{\right}
\def\bg{\bigg}
\def\({\bg(}
\def\[{\bg\lfloor}
\def\){\bg)}
\def\]{\bg\rfloor}
\def\t{\text}
\def\f{\frac}

\def\sm{\setminus}

\def\bi{\binom}
\def\eq{\equiv}

\def\ls{\leqslant}
\def\gs{\geqslant}
\def\mo{\roman{mod}}

\def\ve{\varepsilon}
\def\al{\alpha}
\def\da{\delta}

\def\Proof{\noindent{\it Proof}}

\def\Remark{\medskip\noindent{\it  Remark}}

\def\Ack{\medskip\noindent {\bf Acknowledgment}}
\hbox {Taiwanese J. Math. 17(2013), no.\,5, 1523-1543.}
\medskip

\topmatter
\title Fibonacci numbers modulo cubes of primes\endtitle
\author Zhi-Wei Sun\endauthor
\leftheadtext{Zhi-Wei Sun} \rightheadtext{Fibonacci numbers modulo cubes of primes}
\affil Department of Mathematics, Nanjing University\\
 Nanjing 210093, People's Republic of China
  \\  zwsun\@nju.edu.cn
  \\ {\tt http://math.nju.edu.cn/$\sim$zwsun}
\endaffil
\abstract Let $p$ be an odd prime. It is well known that $F_{p-(\f p5)}\eq0\ (\mo\ p)$, where
$\{F_n\}_{n\gs0}$ is the Fibonacci sequence and $(-)$ is the Jacobi symbol.
In this paper we show that if $p\not=5$ then we may determine $F_{p-(\f p5)}$ mod $p^3$ in the following way:
$$\sum_{k=0}^{(p-1)/2}\f{\bi{2k}k}{(-16)^k}\eq\l(\f{p}5\r)\l(1+\f{F_{p-(\f {p}5)}}2\r)\ \ (\mo\ p^3).$$
We also use Lucas quotients to determine
$\sum_{k=0}^{(p-1)/2}\bi{2k}k/m^k$ modulo $p^2$ for any integer $m\not\eq0\ (\mo\ p)$;
in particular, we obtain
$$\sum_{k=0}^{(p-1)/2}\f{\bi{2k}k}{16^k}\eq\l(\f3{p}\r)\ \ (\mo\ p^2).$$
In addition, we pose three conjectures for further research.
\endabstract
\thanks 2010 {\it Mathematics Subject Classification}.\,Primary 11B39, 11B65;
Secondary 05A10, 11A07.
\newline\indent Supported by the National Natural Science
Foundation (grant 11171140) of China and the PAPD of Jiangsu Higher
Education Institutions.
\endthanks
\endtopmatter
\document

\heading{1. Introduction}\endheading

The well known Fibonacci sequence $\{F_n\}_{n\gs0}$, defined by
$$F_0=0,\ F_1=1,\ \t{and}\ F_{n+1}=F_n+F_{n-1}\ (n=1,2,3,\ldots),$$
plays an important role in many fields of mathematics. This sequence has nice number-theoretic properties;
for example, E. Lucas showed that $(F_m,F_n)=F_{(m,n)}$ for any $m,n\in\N=\{0,1,\ldots\}$, where $(m,n)$ denotes the
greatest common divisor of $m$ and $n$.

Let $p\not=2,5$ be a prime. It is known that $F_{p-(\f p5)}\eq0\ (\mo\ p)$, where $(-)$ denotes
the Jacobi symbol. In 1992 Z. H. Sun and Z. W. Sun [SS] proved that if $p^2\nmid F_{p-(\f p5)}$ then the Fermat equation $x^p+y^p=z^p$
has no integral solutions with $p\nmid xyz$. When $F_{p-(\f p5)}\eq0\ (\mo\ p^2)$, $p$ is called
a Wall-Sun-Sun prime (cf. [CDP] and [CP, p.\,32]).
It is conjectured that there should be infinitely many (but rare) Wall-Sun-Sun primes though none of them has been found.
There are some congruences for the Fibonacci quotient $F_{p-(\f p5)}/p$ modulo $p$ (cf. [W], [SS] and [ST]);
for example, in 1982 H. C. Williams [W] proved that
$$\f{F_{p-(\f p5)}}p\eq\f 25\sum_{k=1}^{\lfloor\f45p\rfloor}\f{(-1)^k}k\pmod{p}.$$
Quite recently H. Pan and Z. W. Sun [PS] proved that for any $a\in\Z^+=\{1,2,3,\ldots\}$ we have
$$\sum_{k=0}^{p^a-1}(-1)^k\bi{2k}k\eq\l(\f{p^a}5\r)\l(1-2F_{p^a-(\f{p^a}5)}\r)\pmod{p^3},$$
which was a conjecture in [ST].

 Now we give the first theorem of this paper.

\proclaim{Theorem 1.1} Let $p$ be an odd prime and let $a$ be a positive integer.
If $p\not=5$, then
$$\sum_{k=0}^{(p^a-1)/2}\f{\bi{2k}k}{(-16)^k}\eq\l(\f{p^a}5\r)\l(1+\f{F_{p^a-(\f{p^a}5)}}2\r)\ (\mo\ p^3).\tag1.1$$
If $p\not=3$, then
$$\sum_{k=0}^{(p^a-1)/2}\f{\bi{2k}k}{(-32)^k}
\eq\l(\f2{p^a}\r)\(1+\f{2^{p^a-1}-1}6-\f{(2^{p^a-1}-1)^2}8\)\ (\mo\ p^3).\tag1.2$$
\endproclaim

Let $p$ be an odd prime and let $a\in\Z^+$. For $k=0,1,\ldots,(p^a-1)/2$, clearly
$$\f{\bi{(p^a-1)/2}k}{\bi{-1/2}k}=\prod_{0\ls j<k}\l(1-\f{p^a}{2j+1}\r)\eq1\pmod{p}$$
and hence
$$\bi{(p^a-1)/2}k\eq\bi{-1/2}k=\f{\bi{2k}k}{(-4)^k}\pmod p.\tag1.3$$
Thus, for any integer $m\not\eq0\ (\mo\ p)$ we have
$$\sum_{k=0}^{p^a-1}\f{\bi{2k}k}{m^k}\eq\sum_{k=0}^{(p^a-1)/2}\f{\bi{2k}k}{m^k}\eq\l(\f{m(m-4)}{p^a}\r)\pmod p,\tag1.4$$
since
$$\sum_{k=0}^{(p^a-1)/2}\bi{(p^a-1)/2}k\l(-\f 4m\r)^k
=\l(1-\f4m\r)^{(p^a-1)/2}$$
and
$$\bi{2k}k=\bi{p^a+(2k-p^a)}{0p^a+k}\eq\bi{2k-p^a}k=0\ \ \ (\mo\ p)$$
for each $k=(p^a+1)/2,\ldots,p^a-1$ by Lucas' theorem (cf. [St, p.\,44]).
Recently the author [Su10] determined
$\sum_{k=0}^{p^a-1}\bi{2k}k/m^k$ mod $p^2$ in terms of Lucas
sequences. See also [SSZ], [GZ] and [Su11a] for related results on
$p$-adic valuations.

Let $A,B\in\Z$. The Lucas sequences $u_n=u_n(A,B)\ (n\in\N)$ and  $v_n=v_n(A,B)\ (n\in\N)$
are defined by
$$u_0=0,\ u_1=1,\ \t{and}\ u_{n+1}=Au_n-Bu_{n-1}\ (n=1,2,3,\ldots)$$
and
$$v_0=2,\ v_1=A,\ \t{and}\ v_{n+1}=Av_n-Bv_{n-1}\ (n=1,2,3,\ldots).$$
The sequence $\{v_n\}_{n\gs0}$ is called the companion of $\{u_n\}_{n\gs0}$.
(Note that $F_n=u_n(1,-1)$, and those $L_n=v_n(1,-1)$ are called Lucas numbers.)
It is known that for any prime $p$ not dividing $2B$ we have
$$u_p\eq\l(\f{\Delta}p\r)\ (\mo\ p)\ \ \ \t{and}\ \ \ u_{p-(\f{\Delta}p)}\eq0\ \ (\mo\ p)$$
where $\Delta=A^2-4B$ (see, e.g., [Su10, Lemma 2.3]); the integer $u_{p-(\f{\Delta}p)}/p$ is called a {\it Lucas quotient}.
The reader may consult [Su06] for connections between Lucas quotients and quadratic fields.

Our second theorem is as follows.

\proclaim{Theorem 1.2} Let $p$ be an odd prime and let $a\in\Z^+$.
Let $m$ be any integer not divisible by $p$.
Then
$$\aligned\sum_{k=0}^{(p^a-1)/2}\f{\bi{2k}k}{m^k}
\eq&\l(\f{m(m-4)}{p^a}\r)
\\&+\l(\f{-m}p\r)\l(\f{m(m-4)}{p^{a-1}}\r)\bar mu_{p-(\f{4-m}p)}(4,m)\ (\mo\ p^2),
\endaligned\tag1.5$$
where
$$\bar m=\cases 1&\t{if}\ m\eq4\ (\mo\ p),
\\2&\t{if}\ (\f{4-m}p)=1,\\2/m&\t{if}\ (\f{4-m}p)=-1.\endcases$$
We also have
$$\sum_{k=0}^{(p^a-1)/2}\f{C_k}{m^k}
\eq\f{4-m}2\sum_{k=0}^{(p^a-1)/2}\f{\bi{2k}k}{m^k}+\f m2-2p\,\da_{a,1}\l(\f{-m}p\r)\ (\mo\ p^2),\tag1.6$$
where $C_k$ denotes the Catalan number $\f1{k+1}\bi{2k}k=\bi{2k}k-\bi{2k}{k+1}$, and the Kronecker symbol $\delta_{s,t}$ takes $1$ or $0$
according as $s=t$ or not.
\endproclaim
\Remark\ 1.1. For any $m\in\Z\sm\{0\}$ and $n\in\N$, the sum $\sum_{k=0}^nk\bi{2k}k/m^k$
is closely related to $\sum_{k=0}^n\bi{2k}k/m^k$ via the identity
$$\sum_{k=0}^n\l(1-\f{m-4}2k\r)\f{\bi{2k}k}{m^k}=(2n+1)\f{\bi{2n}n}{m^n}$$
which can be easily proved by induction.
\medskip

Now we present two consequences of Theorem 1.2.

\proclaim{Corollary 1.1} Let $p$ be an odd prime and let $a\in\Z^+$. Then
$$\sum_{k=0}^{(p^a-1)/2}\f{\bi{2k}k}{8^k}\eq\l(\f2{p^a}\r)\ (\mo\ p^2)\tag1.7$$
and
$$\sum_{k=0}^{(p^a-1)/2}\f{\bi{2k}k}{16^k}\eq\l(\f3{p^a}\r)\ (\mo\ p^2).\tag1.8$$
\endproclaim

\proclaim{Corollary 1.2} Let $p>3$ be a prime. Then
$$\sum_{k=0}^{(p-1)/2}\f{\bi{2k}k}{(2k-1)^2 16^k}\eq\l(\f{-1}p\r)\f{3(\f p3)+1}4\ (\mo\ p^2),\tag1.9$$
that is,
$$\sum_{k=0}^{(p-1)/2}\f{\bi{2k}k}{(2k-1)^216^k}\eq\cases1\ (\mo\ p^2)&\t{if}\ p\eq1\ (\mo\ 12),
\\-1/2\ (\mo\ p^2)&\t{if}\ p\eq5\ (\mo\ 12),\\-1\ (\mo\ p^2)&\t{if}\ p\eq7\ (\mo\ 12),
\\1/2\ (\mo\ p^2)&\t{if}\ p\eq11\ (\mo\ 12).\endcases\tag1.10$$
\endproclaim

We will show Theorems 1.1 and 1.2 in Sections 2 and 3 respectively.
Section 4 is devoted to the proofs of Corollaries 1.1--1.2.
\smallskip

To conclude this section we pose three conjectures.

\proclaim{Conjecture 1.1} For any $n\in\N$ we have
$$\f1{(2n+1)^2\bi{2n}n}\sum_{k=0}^n\f{\bi{2k}k}{16^k}\eq\cases1\ (\mo\ 9)&\t{if}\ 3\mid n,\\4\ (\mo\ 9)&\t{if}\ 3\nmid n.
\endcases$$ Also,
$$\f1{3^{2a}}\sum_{k=0}^{(3^a-1)/2}\f{\bi{2k}k}{16^k}\eq(-1)^a10\ (\mo\ 27)$$
for every $a=1,2,3,\ldots.$
\endproclaim

Let $p>3$ be a prime. In 2007 A. Adamchuk [A] conjectured that if $p\eq1\ (\mo\ 3)$ then
$$\sum_{k=1}^{\lfloor\f23p\rfloor}\bi{2k}k\eq0\ (\mo\ p^2).$$
Motivated by this and Theorems 1.1 and 1.2, we pose the following conjecture based on the author's computation
via the software {\tt Mathematica}.

\proclaim{Conjecture 1.2} Let $p$ be an odd prime and let $a\in\Z^+$.

{\rm (i)} If $p\eq1\ (\mo\ 3)$ or $a>1$, then
$$\sum_{k=0}^{\lfloor\f56p^a\rfloor}\f{\bi{2k}k}{16^k}\eq\l(\f 3{p^a}\r)\ (\mo\ p^2).$$

{\rm (ii)} Suppose $p\not=5$. If $p^a\eq1,2\ (\mo\ 5)$ or $p\eq2\ (\mo\ 5)$ or $a>2$, then
$$\sum_{k=0}^{\lfloor\f 45p^a\rfloor}(-1)^k\bi{2k}k\eq\l(\f5{p^a}\r)\ (\mo\ p^2).$$
If $p^a\eq1,3\ (\mo\ 5)$ or $p\eq3\ (\mo\ 5)$ or $a>2$, then
$$\sum_{k=0}^{\lfloor\f 35p^a\rfloor}(-1)^k\bi{2k}k\eq\l(\f5{p^a}\r)\ (\mo\ p^2).$$

{\rm (iii)} If $p^a\eq1,2\ (\mo\ 5)$ or $p\eq2\ (\mo\ 5)$ or $a>2$, then
$$\sum_{k=0}^{\lfloor\f 7{10}p^a\rfloor}\f{\bi{2k}k}{(-16)^k}\eq\l(\f 5{p^a}\r)\ (\mo\ p^2).$$
If $p^a\eq1,3\ (\mo\ 5)$ or $p\eq3\ (\mo\ 5)$ or $a>2$, then
$$\sum_{k=0}^{\lfloor\f 9{10}p^a\rfloor}\f{\bi{2k}k}{(-16)^k}\eq\l(\f 5{p^a}\r)\ (\mo\ p^2).$$
\endproclaim

\proclaim{Conjecture 1.3} Let $p\not=2,5$ be a prime and set $q:=F_{p-(\f p5)}/p$. Then
$$p\sum_{k=1}^{p-1}\f{F_{2k}}{k^2\bi{2k}k}\eq-\l(\f p5\r)\l(\f 32q+\f 54p\, q^2\r)\ \ (\mo\ p^2)$$
and $$p\sum_{k=1}^{p-1}\f{L_{2k}}{k^2\bi{2k}k}\eq-\f 52q-\f {15}4p\, q^2\ \ \ (\mo\ p^2).$$
\endproclaim
\Remark\ 1.2. It is interesting to compare Conjecture 1.3 with the two identities
$$\sum_{k=1}^\infty\f{F_{2k}}{k^2\bi{2k}k}=\f{4\pi^2}{25\sqrt 5}\ \ \t{and}\ \ \sum_{k=1}^\infty\f{L_{2k}}{k^2\bi{2k}k}=\f{\pi^2}5$$
obtained by putting $x=(\sqrt5\pm1)/2$ in the known formula
$$2\arcsin^2\f x2=\sum_{k=1}^\infty\f{x^{2k}}{k^2\bi{2k}k}\ \ \ (|x|<2).$$

\heading{2. Proof of Theorem 1.1}\endheading

\proclaim{Lemma 2.1} Let $p$ be an odd prime and let $k\in\{0,\ldots,(p^a-1)/2\}$ with $a\in\Z^+$. Then
$$\aligned&\bi{(p^a-1)/2+k}{2k}-\f{\bi{2k}k}{(-16)^k}
\\\eq&(-1)^{k-1}\l(\f{-1}{p^a}\r)\bi{p^a-1-2k}{(p^a-1)/2-k}\sum_{0<j\ls k}\f{p^{2a}}{(2j-1)^2}\pmod{p^3}.
\endaligned
\tag2.1$$
\endproclaim
\Proof. Clearly (2.1) holds for $k=0$. Below we assume $1\ls k\ls (p^a-1)/2$.
Note that
$$\align\bi{(p^a-1)/2+k}{2k}=&\f{\prod_{j=1}^k(p^{2a}-(2j-1)^2)}{4^k(2k)!}
\\=&\f{\prod_{j=1}^k(-(2j-1)^2)}{4^k(2k)!}\prod_{j=1}^k\l(1-\f{p^{2a}}{(2j-1)^2}\r)
\\\eq&\f{\bi{2k}k}{(-16)^k}\(1-\sum_{j=1}^k\f{p^{2a}}{(2j-1)^2}\)\ \ (\mo\ p^4)
\endalign$$
(which was observed by Z. H. Sun [S11, Lemma 2.2] in the case $a=1$).
Thus, in view of (1.3) we have
$$\align\f{\bi{2k}k}{(-16)^k}
\eq&\bi{(p^a-1)/2}k4^{-k}=\bi{(p^a-1)/2}{(p^a-1)/2-k}4^{-k}
\\\eq&\bi{p^a-1-2k}{(p^a-1)/2-k}(-4)^{-((p^a-1)/2-k)}4^{-k}
\\\eq&\l(\f{-1}{p^a}\r)(-1)^k\bi{p^a-1-2k}{(p^a-1)/2-k}\pmod{p}.
\endalign$$
So (2.1) holds. \qed

\proclaim{Lemma 2.2} Let $p$ be an odd prime and let $a\in\Z^+$. If $p\not=3$, then
$$1+\f{2^{p^a-1}-1}6+\f{(2^{p^a-1}-1)^2}{24}\eq\l(\f2{p^a}\r)\f{2^{p^a}+1}{3\times2^{(p^a-1)/2}}\pmod{p^3}.\tag2.2$$
When $p\not=5$, we have
$$\f{L_{p^a}-1}5-\l(\f{p^a}5\r)F_{p^a}+1\eq-\f12F_{p^a-(\f{p^a}5)}^2\pmod{p^4}.\tag2.3$$
\endproclaim
\Proof. Note that
$$2^{(p^a-1)/2}=\l(2^{\f{p-1}2}\r)^{\sum_{k=0}^{a-1}p^k}\eq\l(\f2p\r)^a=\l(\f{2}{p^a}\r)\pmod{p}$$
and
$$\align\f{2^{p^a-1}-1}p=&\f{2^{(p^a-1)/2}-(\f2{p^a})}p\l(2^{(p^a-1)/2}+\l(\f2{p^a}\r)\r)
\\\eq& 2\l(\f 2{p^a}\r)\f{2^{(p^a-1)/2}-(\f2{p^a})}p\pmod{p}.
\endalign$$
Thus
$$\align(2^{p^a-1}-1)^2\eq&4\l(2^{(p^a-1)/2}-\l(\f2{p^a}\r)\r)^2
\\=&4(2^{p^a-1}-1)+8-8\l(\f2{p^a}\r)2^{(p^a-1)/2}\pmod{p^3}.
\endalign$$
It follows that
$$\align&\f{2^{(p^a-1)/2}}2\l(2^{p^a-1}-1\r)+\f18\l(\f2{p^a}\r)\l(2^{p^a-1}-1\r)^2
\\\eq&\f{2^{(p^a-1)/2}}8\l(\l(2^{p^a-1}-1\r)^2-8+8\l(\f 2{p^a}\r)2^{(p^a-1)/2}\r)
\\&+\f18\l(\f2{p^a}\r)\l(2^{p^a-1}-1\r)^2
\\\eq&\f14\l(\f2{p^a}\r)\l(2^{p^a-1}-1\r)^2-2^{(p^a-1)/2}+\l(\f2{p^a}\r)2^{p^a-1}
\\\eq&\l(\f2{p^a}\r)\l(2^{p^a}+1\r)-3\times2^{(p^a-1)/2}\ \ \pmod{p^3},
\endalign$$
which is equivalent to (2.2) times $3\times2^{(p^a-1)/2}$. So (2.2) is valid if $p>3$.

Now assume that $p\not=5$.
As $L_{2n}=5F_n^2+2(-1)^n=L_n^2-2(-1)^n$ for all $n\in\N$, by [SS, Corollary 1] we have
$L_{p-(\f 5p)}\eq 2(\f p5)\pmod{p^2}$. Thus, in view of [Su10, Lemma 2.3],
$$L_{p^a-(\f{p^a}5)}\eq(-1)^{((\f5p)-(\f 5{p^a}))/2}L_{p-(\f 5p)}
\eq 2\l(\f{p^a}5\r)\pmod{p^2}.$$
Write $L_{p^a-(\f{p^a}5)}=2(\f{p^a}5)+p^2Q$ with $Q\in\Z$. Then
$$\align 5F_{p^a-(\f{p^a}5)}^2=&L_{p^a-(\f{p^a}5)}^2-4(-1)^{p^a-(\f{p^a}5)}
\\=&-4+\l(2\l(\f{p^a}5\r)+p^2Q\r)^2\eq 4p^2\l(\f {p^a}5\r)Q\ \pmod{p^4}.
\endalign$$
Note that
$$L_{p^a}=F_{p^a}+2F_{p^a-1}=2F_{p^a+1}-F_{p^a}=2F_{p^a-(\f {p^a}5)}+\l(\f{p^a}5\r)F_{p^a}$$
and
$$2L_{p^a}=5F_{p^a-1}+L_{p^a-1}=5F_{p^a+1}-L_{p^a+1}=5F_{p^a-(\f{p^a}5)}+\l(\f{p^a}5\r)L_{p^a-(\f{p^a}5)}.$$
Therefore
$$\align&\f{L_{p^a}-1}5-\l(\f{p^a}5\r)F_{p^a}+1
\\=&2F_{p^a-(\f{p^a}5)}-\f45(L_{p^a}-1)
\\=&2F_{p^a-(\f{p^a}5)}-\f45\l(\f52F_{p^a-(\f{p^a}5)}+\f12\l(\f{p^a}5\r)L_{p^a-(\f{p^a}5)}-1\r)
\\=&\f 45-\f25\l(\f{p^a}5\r)\l(2\l(\f{p^a}5\r)+p^2Q\r)\eq-\f12F_{p^a-(\f{p^a}5)}^2\pmod{p^4}.
\endalign$$
This proves (2.3). \qed
\smallskip

The following Lemma was posed as [Su12, Conjecture 1.1].

 \proclaim{Lemma 2.3} Let $p$ be an odd prime and let $H_k^{(2)}=\sum_{0<j\ls k}1/j^2$ for $k=0,1,2,\ldots$.
If $p\not=5$, then
$$\sum_{k=0}^{p-1}(-1)^k\bi{2k}kH_k^{(2)}\eq\l(\f p5\r)\f52q^2-\delta_{p,3}\pmod{p},\tag2.4$$
where $q$ denotes the Fibonacci quotient $F_{p-(\f p5)}/p$.
If $p>3$, then
$$\sum_{k=0}^{p-1}(-2)^k\bi{2k}kH_k^{(2)}\eq\f23q_p(2)^2\pmod{p},\tag2.5$$
where $q_p(2)$ denotes the Fermat quotient $(2^{p-1}-1)/p$.
\endproclaim
\Proof. It is easy to verify (2.4) for $p=3$. Below we assume $p>3$.

The desired congruences essentially follow from [MT, (37)]. Here we provide the details.
Putting $t=-1,-1/2$ in [MT, (37)] we get
$$\sum_{k=0}^{p-1}(-1)^k\bi{2k}kH_k^{(2)}\eq-2\sum_{k=1}^{p-1}\f{u_k(3,1)}{k^2}\pmod p\tag2.6$$
and
$$\sum_{k=0}^{p-1}(-2)^k\bi{2k}kH_k^{(2)}\eq-2\sum_{k=1}^{p-1}\f{u_k(5/2,1)}{k^2}\pmod p.\tag2.7$$
Note that $u_k(3,1)=F_{2k}$ and $u_k(5/2,1)=2(2^k-2^{-k})/3$
for all $k=0,1,2,\ldots$.

Since
$$\sum_{k=1}^{p-1}\f{2^k}{k^2}\eq-q_p(2)^2\pmod p\ \ \t{and}\ \ \sum_{k=1}^{p-1}\f1{k^22^k}\eq-\f{q_p(2)^2}2\pmod p$$
by [Gr] and [S08, Theorem 4.1(iv)] respectively,  (2.7) implies that
$$\sum_{k=0}^{p-1}(-2)^k\bi{2k}kH_k^{(2)}\eq-\f43\sum_{k=1}^{p-1}\l(\f{2^k}{k^2}-\f{1}{k^22^k}\r)\eq\f 23q_p(2)^2\pmod p.$$

Now we work with $p>5$. Recall that for any $n\in\N$ we have
$$F_n=\f{\al^n-\beta^n}{\al-\beta}\quad\t{and}\quad L_n=\al^n+\beta^n,$$
where $\al$ and $\beta$ are the two roots of the equation $x^2-x-1=0$.
By [PS, (3.2) and (3.7)],
$$\align 2\beta^{2p}\sum_{k=1}^{p-1}\f{\al^{2k}}{k^2}\eq&-2\sum_{k=1}^{p-1}\f{\al^k}{k^2}-\l(\f{L_p-1}p\r)^2
\\\eq&\l(\f{2\al^p-1}5-1\r)\l(\f{L_p-1}p\r)^2\pmod p.
\endalign$$
Since $\al\beta=-1$ and $\al^{2p}=(\al+1)^p\eq\al^p+1\ (\mo\ p)$, we have
$$\sum_{k=1}^{p-1}\f{\al^{2k}}{k^2}\eq\f{(\al^p+1)(\al^p-3)}5\l(\f{L_p-1}p\r)^2\eq-\f{\al^p+2}5\l(\f{L_p-1}p\r)^2\pmod{p}.$$
Hence
$$\align 5\sum_{k=1}^{p-1}\f{F_{2k}}{k^2}=&(\al-\beta)^2\sum_{k=1}^{p-1}\f{F_{2k}}{k^2}=(\al-\beta)\sum_{k=1}^{p-1}\f{\al^{2k}-\beta^{2k}}{k^2}
\\\eq&(\al-\beta)\f{\beta^p-\al^p}5\l(\f{L_p-1}p\r)^2\eq-\f{(\al-\beta)^{p+1}}5\l(\f{L_p-1}p\r)^2
\\=&-5^{(p-1)/2}\l(\f{L_p-1}p\r)^2\eq-\l(\f 5p\r)\f{25}4q^2\pmod{p}
\endalign$$
since $2(L_p-1)\eq 5F_{p-(\f p5)}\ (\mo\ p^2)$ by [ST, p.\,139].
Combining this with (2.6) we obtain
$$\sum_{k=0}^{p-1}(-1)^k\bi{2k}kH_k^{(2)}\eq2\times\l(\f 5p\r)\f{5}4q^2=\l(\f p5\r)\f 52q^2\pmod{p}.$$

The proof of Lemma 2.3 is now complete. \qed

\medskip
\noindent{\it Proof of Theorem 1.1}. Let us first recall the following two identities:
$$F_{n+1}=\sum_{k=0}^{\lfloor n/2\rfloor}\bi{n-k}k\ \ \ \t{and}\ \ \
\sum_{k=0}^n\f{\bi{n+k}{2k}}{2^k}=\f{2^{2n+1}+1}{3\times2^n}.$$
Thus we have
$$\align F_{p^a}=&\sum_{k=0}^{(p^a-1)/2}\bi{p^a-1-k}{p^a-1-2k}=\sum_{j=0}^{(p^a-1)/2}\bi{(p^a-1)/2+j}{2j}
\endalign$$
and
$$\sum_{k=0}^{(p^a-1)/2}\f{\bi{(p^a-1)/2+k}{2k}}{2^k}=\f{2^{p^a}+1}{3\times2^{(p^a-1)/2}}.$$
Therefore, with the help of (2.1),
$$\align&\sum_{k=0}^{(p^a-1)/2}\f{\bi{2k}k}{(-16)^k}-F_{p^a}
\\\eq&\sum_{k=0}^{(p^a-1)/2}\bi{p^a-1-2k}{(p^a-1)/2-k}(-1)^{(p^a-1)/2-k}\sum_{0<j\ls k}\f{p^{2a}}{(2j-1)^2}
\\=&\sum_{k=0}^{(p^a-1)/2}\bi{2k}k(-1)^k\sum_{0<j\ls (p^a-1)/2-k}\f{p^{2a}}{(2j-1)^2}
\pmod{p^3}\endalign$$
and
$$\align&\sum_{k=0}^{(p^a-1)/2}\f{\bi{2k}k}{(-32)^k}-\f{2^{p^a}+1}{3\times2^{(p^a-1)/2}}
\\\eq&\sum_{k=0}^{(p^a-1)/2}\bi{2k}k\f{(-1)^k}{2^{(p^a-1)/2-k}}\sum_{0<j\ls(p^a-1)/2-k}\f{p^{2a}}{(2j-1)^2}
\pmod{p^3}.\endalign$$

For $k=0,\ldots,(p^a-1)/2$, clearly
$$\align\sum_{j=1}^{(p^a-1)/2-k}\f{p^{2a}}{(2j-1)^2}=&
\sum_{j=k+1}^{(p^a-1)/2}\f{p^{2a}}{(2((p^a-1)/2-j+1)-1)^2}
\\\eq&\sum_{j=1}^{(p^a-1)/2}\f{p^{2a}}{4j^2}-\sum_{0<j\ls k}\f{p^{2a}}{4j^2}
\\\eq&\sum_{i=1}^{(p-1)/2}\f{p^{2a}}{4(p^{a-1}i)^2}
-\sum_{0<i\ls\lfloor k/p^{a-1}\rfloor}\f{p^{2a}}{4(p^{a-1}i)^2}\ (\mo\ p^3).
\endalign$$
Since
$$2\sum_{i=1}^{(p-1)/2}\f1{i^2}\eq\sum_{i=1}^{(p-1)/2}\l(\f1{i^2}+\f1{(p-i)^2}\r)=\sum_{i=1}^{p-1}\f1{i^2}\eq2\delta_{p,3}\ \ (\mo\ p)$$
with the help of Wolstenholme's congruence (cf. [Wo] and [Z]),
by the above we have
$$\aligned&\f{-4}{p^2}\(\sum_{k=0}^{(p^a-1)/2}\f{\bi{2k}k}{(-16)^k}-F_{p^a}\)
\\\eq&\sum_{k=0}^{p^a-1}\bi{2k}k(-1)^k\l(H_{\lfloor k/p^{a-1}\rfloor}^{(2)}-\delta_{p,3}\r)\ \ \ (\mo\ p)
\endaligned\tag2.8$$
and
$$\aligned&\f{-4}{p^2}\(\sum_{k=0}^{(p^a-1)/2}\f{\bi{2k}k}{(-32)^k}-\f{2^{p^a}+1}{3\times2^{(p^a-1)/2}}\)
\\\eq&\l(\f2{p^a}\r) \sum_{k=0}^{p^a-1}\bi{2k}k(-2)^kH_{\lfloor k/p^{a-1}\rfloor}^{(2)}\ \ \pmod{p}.
\endaligned\tag2.9$$

Our next task is to simplify the right-hand sides of congruences (2.8) and (2.9).
Let $u\in\{1,2\}$. Then
$$\align&\sum_{k=0}^{p^a-1}\bi{2k}k(-u)^kH_{\lfloor k/p^{a-1}\rfloor}^{(2)}
\\=&\sum_{k=0}^{p-1}\sum_{r=0}^{p^{a-1}-1}\bi{2p^{a-1}k+2r}{p^{a-1}k+r}(-u)^{p^{a-1}k+r}H_k^{(2)}
\\\eq&\sum_{k=0}^{p-1}(-u)^kH_k^{(2)}\sum_{r=0}^{p^{a-1}-1}\bi{2p^{a-1}k+2r}{p^{a-1}k+r}(-u)^r
\ \pmod{p}.
\endalign$$
For $k\in\{0,\ldots,p-1\}$ and $r\in\{0,\ldots,p^{a-1}-1\}$, by the Chu-Vandermonde identity (cf. [GKP, p.\,169]) we have
$$\bi{2p^{a-1}k+2r}{p^{a-1}k+r}=\sum_{j=0}^{p^{a-1}k+r}\bi{2p^{a-1}k}{j}\bi{2r}{p^{a-1}k+r-j}.$$
If $p^{a-1}\nmid j$, then
$$\bi{2p^{a-1}k}j=\f{2p^{a-1}k}j\bi{2p^{a-1}k-1}{j-1}\eq0\ \ (\mo\ p).$$
Thus
$$\align\bi{2p^{a-1}k+2r}{p^{a-1}k+r}
\eq&\sum_{j=0}^k\bi{2p^{a-1}k}{p^{a-1}j}\bi{2r}{p^{a-1}(k-j)+r}=\bi{2p^{a-1}k}{p^{a-1}k}\bi{2r}r
\\\eq&\bi{2k}k\bi{2r}r\ \pmod{p}\quad\t{(by Lucas' theorem)}.\endalign$$
Therefore
$$\aligned&\sum_{k=0}^{p^a-1}\bi{2k}k(-u)^kH_{\lfloor k/p^{a-1}\rfloor}^{(2)}
\\\eq&\sum_{k=0}^{p-1}(-u)^kH_k^{(2)}\bi{2k}k\sum_{r=0}^{p^{a-1}-1}\bi{2r}r(-u)^r
\ \pmod{p}.
\endaligned\tag2.10$$

 In view of (1.4),
$$\sum_{r=0}^{p^{a-1}-1}\bi{2r}r(-1)^r\eq\l(\f5{p^{a-1}}\r)\pmod p,$$
and also
$$\sum_{r=0}^{p^{a-1}-1}\bi{2r}r(-2)^r\eq\l(\f{(p-1)/2\times((p-1)/2-4)}{p^{a-1}}\r)=1\pmod{p}$$
provided $p\not=3$. Combining this with (2.8) and (2.10), we obtain
$$\aligned&\f{-4}{p^2}\(\sum_{k=0}^{(p^a-1)/2}\f{\bi{2k}k}{(-16)^k}-F_{p^a}\)
\\\eq&\sum_{k=0}^{p-1}\bi{2k}k(-1)^kH_k^{(2)}\sum_{r=0}^{p^{a-1}-1}\bi{2r}r(-1)^r-\delta_{p,3}\sum_{k=0}^{p^a-1}\bi{2k}k(-1)^k
\\\eq&\sum_{k=0}^{p-1}\bi{2k}k(-1)^kH_k^{(2)}\l(\f{5}{p^{a-1}}\r)-\delta_{p,3}\l(\f 5{p^a}\r)\pmod{p},
\endaligned$$
and hence
$$\aligned&\f{-4}{p^2}\(\sum_{k=0}^{(p^a-1)/2}\f{\bi{2k}k}{(-16)^k}-F_{p^a}\)
\\\eq&\l(\f{5}{p^{a-1}}\r)\(\sum_{k=0}^{p-1}\bi{2k}k(-1)^kH_k^{(2)}+\da_{p,3}\)\pmod{p}.
\endaligned\tag2.11$$
Similarly, when $p\not=3$ we have
$$\aligned&\f{-4}{p^2}\(\sum_{k=0}^{(p^a-1)/2}\f{\bi{2k}k}{(-32)^k}-\f{2^{p^a}+1}{3\times2^{(p^a-1)/2}}\)
\\\eq&\l(\f2{p^a}\r)\sum_{k=0}^{p-1}\bi{2k}k(-2)^kH_k^{(2)}\pmod{p}.
\endaligned\tag2.12$$

Now assume that $p\not=3$. By (2.5) and (2.12),
$$\sum_{k=0}^{(p^a-1)/2}\f{\bi{2k}k}{(-32)^k}-\f{2^{p^a}+1}{3\times2^{(p^a-1)/2}}
\eq-\l(\f2{p^a}\r)\f{p^2}6q_p(2)^2\pmod{p^3}.$$
Since $p^a\eq p\pmod{\varphi(p^2)}$, we have $2^{p^a}\eq 2^p\pmod{p^2}$ and hence
$$\sum_{k=0}^{(p^a-1)/2}\f{\bi{2k}k}{(-32)^k}-\f{2^{p^a}+1}{3\times2^{(p^a-1)/2}}
\eq-\l(\f2{p^a}\r)\f{(2^{p^a-1}-1)^2}6\pmod{p^3}.$$
Combining this with (2.2) we immediately obtain (1.2).

Below we suppose that $p\not=5$. By (2.4) and (2.11),
$$\sum_{k=0}^{(p^a-1)/2}\f{\bi{2k}k}{(-16)^k}-F_{p^a}\eq-\f{5}8\l(\f{p^a}5\r)F_{p-(\f p5)}^2\pmod{p^3}.$$
In view of [Su10, Lemma 2.3],
$$\f{F_{p^a-(\f{p^a}5)}}p\eq (-1)^{((\f5p)-(\f5{p^a}))/2}\l(\f 5{p^{a-1}}\r)
\f{F_{p-(\f p5)}}p=\f{F_{p-(\f p5)}}p\ (\mo\ p)$$
and thus
$$\sum_{k=0}^{(p^a-1)/2}\f{\bi{2k}k}{(-16)^k}-F_{p^a}\eq-\f{5}8\l(\f{p^a}5\r)F_{p^a-(\f {p^a}5)}^2\pmod{p^3}.$$
Combining this with (2.3) we get
$$\sum_{k=0}^{(p^a-1)/2}\f{\bi{2k}k}{(-16)^k}-F_{p^a}\eq\f{L_{p^a}-1}4\l(\f{p^a}5\r)-\f 54F_{p^a}+\f 54\l(\f{p^a}5\r)
\pmod{p^3}.$$
Therefore
$$\align \sum_{k=0}^{(p^a-1)/2}\f{\bi{2k}k}{(-16)^k}\eq&\l(\f{p^a}5\r)\f{L_{p^a}}4-\f{F_{p^a}}4+\l(\f{p^a}5\r)
\\=&\l(\f{p^a}5\r)\l(1+\f14\l(L_{p^a}-\l(\f{p^a}5\r)F_{p^a}\r)\r)
\\=&\l(\f{p^a}5\r)\l(1+\f12F_{p^a-(\f{p^a}5)}\r)
\pmod{p^3}.
\endalign$$
This proves (1.1).

So far we have completed the proof of Theorem 1.1. \qed

\heading{3. Proof of Theorem 1.2}\endheading

  We need some preliminary results about Lucas sequences.

 Let $A,B\in\Z$ and $\Delta=A^2-4B$.
The equation $x^2-Ax+B=0$ has two roots
$$\al=\f{A+\sqrt{\Delta}}2\quad\t{and}\quad\beta=\f{A-\sqrt{\Delta}}2$$
which are algebraic integers.
It is well known that for any $n\in\N$ we have
$$u_n(A,B)=\sum_{0\ls k<n}\al^k\beta^{n-1-k}\quad\t{and }\quad v_n(A,B)=\al^n+\beta^n.$$
If $p$ is a prime then
$$v_p(A,B)=\al^p+\beta^p\eq(\al+\beta)^p=A^p\eq A\ (\mo\ p).$$

\proclaim{Lemma 3.1} Let $A,B\in\Z$ and $n\in\N$. Then
$$u_{n+1}(A,B)=\sum_{k=0}^{\lfloor n/2\rfloor}\bi{n-k}kA^{n-2k}(-B)^k.\tag3.1$$
\endproclaim
\Remark\ 3.1. (3.1) is a well-known formula due to Lagrange,  see, e.g., H. Gould [G, (1.60)].

\proclaim{Lemma 3.2} Let $A,B\in\Z$ and let $p$ be an odd prime not dividing $B\Delta$ where $\Delta=A^2-4B$.
Then
$$u_p(A,B)\eq\f A2B^{((\f{\Delta}p)-1)/2}u_{p-(\f{\Delta}p)}(A,B)+\l(\f{\Delta}p\r)\f{B^{p-1}+1}2\ (\mo\ p^2).\tag3.2$$
\endproclaim
\Proof. For convenience we let $u_n=u_n(A,B)$ and $v_n=v_n(A,B)$ for all $n\in\N$.

Let $\al$ and $\beta$ be the two roots of the equation $x^2-Ax+B=0$. Then
$$v_n^2-\Delta u_n^2=(\al^n+\beta^n)^2-\l(\al^n-\beta^n\r)^2=4(\al\beta)^n=4B^n$$
for any $n\in\N$.  As $p\mid u_{p-(\f{\Delta}p)}$ (see, e.g., [Su10, Lemma 2.3]),
$p^2$ divides
$$\align &v_{p-(\f{\Delta}p)}^2-4B^{p-(\f{\Delta}p)}
\\=&\l(v_{p-(\f{\Delta}p)}-2\l(\f Bp\r)B^{(p-(\f{\Delta}p))/2}\r)
\l(v_{p-(\f{\Delta}p)}+2\l(\f Bp\r)B^{(p-(\f{\Delta}p))/2}\r).
\endalign$$
On the other hand, by [Su10, Lemma 2.3] we have
$$v_{p-(\f{\Delta}p)}\eq2B^{(1-(\f{\Delta}p))/2}\eq2\l(\f Bp\r)B^{(p-(\f{\Delta}p))/2}\ (\mo\ p).$$
Therefore
$$v_{p-(\f{\Delta}p)}\eq2\l(\f Bp\r)B^{(p-(\f{\Delta}p))/2}\ (\mo\ p^2).$$

By induction, for $\ve=\pm1$ we have
$$Au_n+\ve v_n=2B^{(1-\ve)/2}u_{n+\ve}$$
for all $n\in\Z^+$. Thus
$$\align 2B^{(1-(\f{\Delta}p))/2}u_p=&Au_{p-(\f{\Delta}p)}+\l(\f{\Delta}p\r)v_{p-(\f{\Delta}p)}
\\\eq&Au_{p-(\f{\Delta}p)}+\l(\f{\Delta}p\r)2\l(\f Bp\r)B^{(p-(\f{\Delta}p))/2}\ \ (\mo\ p^2)
\endalign$$
and hence
$$\align&2u_p-AB^{((\f{\Delta}p)-1)/2}u_{p-(\f{\Delta}p)}
\\\eq&\l(\f{\Delta}p\r)\(2\l(\f Bp\r)\l(B^{(p-1)/2}-\l(\f Bp\r)\r)+2\)
\\\eq&\l(\f{\Delta}p\r)\l(B^{p-1}-1+2\r)\ \ (\mo\ p^2).
\endalign$$
So (3.2) is valid. \qed

\proclaim{Lemma 3.3} Let $p$ be an odd prime and let $a\in\Z^+$. Let $m$ be an integer not divisible by $p$. Then
$$\sum_{k=0}^{(p^a-1)/2}\f{\bi{2k}{k+1}}{m^k}\eq\f{m-2}2\sum_{k=0}^{(p^a-1)/2}\f{\bi{2k}k}{m^k}
-\f m2+2p\da_{a,1}\l(\f{-m}{p}\r)\ (\mo\ p^2).\tag3.3$$
\endproclaim
\Proof. Observe that
$$\align &\sum_{k=0}^{(p^a-1)/2}\f{\bi{2k}k+\bi{2k}{k+1}}{m^k}
\\=&\sum_{k=0}^{(p^a-1)/2}\f{\bi{2k+1}{k+1}}{m^k}
=\f{\bi {p^a}{(p^a+1)/2}}{m^{(p^a-1)/2}}+\f12\sum_{k=0}^{(p^a-3)/2}\f{\bi{2k+2}{k+1}}{m^k}
\\=&\f {p^a/m^{(p^a-1)/2}}{(p^a+1)/2}\bi{p^a-1}{(p^a-1)/2}+\f m2\sum_{k=1}^{(p^a-1)/2}\f{\bi{2k}{k}}{m^k}
\\\eq&2p\da_{a,1}\l(\f{-m}p\r)+\f m2\sum_{k=0}^{(p^a-1)/2}\f{\bi{2k}k}{m^k}-\f m2\ (\mo\ p^2)
\endalign$$
and hence (3.3) follows. \qed

\medskip
\noindent{\it Proof of Theorem 1.2}. Set $n=(p^a-1)/2$. By Lemma 3.3,
$$\sum_{k=0}^n\f{C_k}{m^k}\eq\l(1-\f{m-2}2\r)\sum_{k=0}^n\f{\bi{2k}k}{m^k}+\f m2-2p\da_{a,1}\l(\f{-m}p\r)\ (\mo\ p^2).$$
This proves (1.6). It remains to show (1.5).

By Lemmas 3.1 and 2.1,
$$\align u_{p^a}(4,m)=&\sum_{k=0}^n\bi{2n-k}k4^{2n-2k}(-m)^k
\\=&\sum_{k=0}^n\bi{2n-k}{2(n-k)}16^{n-k}(-m)^k
=\sum_{k=0}^n\bi{n+k}{2k}16^k(-m)^{n-k}
\\\eq&\sum_{k=0}^n\bi{2k}k(-1)^k(-m)^{n-k}=(-m)^n\sum_{k=0}^{n}\f{\bi{2k}k}{m^k}\ (\mo\ p^2).
\endalign$$
Note that
$$(-m)^n=\l((-m)^{(p-1)/2}\r)^{\sum_{s=0}^{a-1}p^s}\eq\l(\f{-m}p\r)^{\sum_{s=0}^{a-1}p^s}=\l(\f{-m}{p^a}\r)\ (\mo\ p)$$
and hence
$$\align(-m)^n-\l(\f{-m}{p^a}\r)\eq&\((-m)^n-\l(\f{-m}{p^a}\r)\)\f{(-m)^n+\l(\f{-m}{p^a}\r)}{2(\f{-m}{p^a})}
\\\eq&\f{(-m)^{p^a-1}-1}2\l(\f{-m}{p^a}\r)\ \ (\mo\ p^2).
\endalign$$
Thus
$$(-m)^n\eq\l(\f{-m}{p^a}\r)\l(1+\f{m^{p^a-1}-1}2\r)\eq\f{(\f{-m}{p^a})}{1-(m^{p^a-1}-1)/2}\ \ (\mo\ p^2)$$
and hence
$$\align\sum_{k=0}^n\f{\bi{2k}k}{m^k}\eq& u_{p^a}(4,m)\l(\f{-m}{p^a}\r)\l(1-\f{m^{p^a-1}-1}2\r)
\\\eq&u_{p^a}(4,m)\l(\f{-m}{p^a}\r)\l(1-\f{m^{p-1}-1}2\r)\ \pmod{p^2}
\endalign$$
since $m^{p^a-1}\eq m^{p-1}\ (\mo\ p^2)$ by Euler's theorem.
By [Su10, Lemma 2.3],
$$\align u_{p^a}(4,m)\eq&\l(\f{4^2-4m}{p^{a-1}}\r)u_p(4,m)\ (\mo\ p^2)
\\\eq&\l(\f{4^2-4m}{p^a}\r)u_1(4,m)=\l(\f{4-m}{p^a}\r)\ (\mo\ p).
\endalign$$
Therefore
$$\align\sum_{k=0}^n\f{\bi{2k}k}{m^k}\eq&u_{p^a}(4,m)\l(\f{-m}{p^a}\r)-u_{p^a}(4,m)\l(\f{-m}{p^a}\r)\f{m^{p-1}-1}2
\\\eq&\l(\f{4-m}{p^{a-1}}\r)\l(\f{-m}{p^a}\r)u_p(4,m)-\l(\f{-m(4-m)}{p^a}\r)\f{m^{p-1}-1}2
\\=&\l(\f{-m}p\r)\l(\f{m(m-4)}{p^{a-1}}\r)u_p(4,m)-\l(\f{m(m-4)}{p^a}\r)\f{m^{p-1}-1}2\ (\mo\ p^2).
\endalign$$

In view of Lemma 3.2,
$$u_p(4,m)-\l(\f{4-m}p\r)\f{m^{p-1}-1}2
\eq\bar mu_{p-(\f{4-m}p)}(4,m)+\l(\f{4-m}p\r)\ (\mo\ p^2).$$
So, by the above, $\sum_{k=0}^n\bi{2k}k/{m^k}$ is congruent to
$$\align &\l(\f{m(m-4)}{p^{a-1}}\r)\l(\f{-m}p\r)
\(\bar mu_{p-(\f{4-m}p)}(4,m)+\l(\f{4-m}p\r)\)
\\=&\l(\f{m(m-4)}{p^a}\r)+\l(\f{-m}p\r)\l(\f{m(m-4)}{p^{a-1}}\r)\bar mu_{p-(\f{4-m}p)}(4,m)
\endalign$$
modulo $p^2$.
This proves (1.5). We are done. \qed

\heading{4. Proofs of Corollaries 1.1--1.2}\endheading

\medskip
\noindent{\it Proof of Corollary 1.1}. Note that $n=p-(\f{4-8}p)\eq0\ (\mo\ 4)$.
The equation $x^2-4x+8=0$ has two roots $2\pm2i$ where $i=\sqrt{-1}$. Thus
$$u_{n}(4,8)=\f{(2+2i)^n-(2-2i)^n}{4i}=\f{(i(2-2i))^n-(2-2i)^n}{4i}=0$$
and hence by Theorem 1.2 we have
$$\sum_{k=0}^{(p^a-1)/2}\f{\bi{2k}k}{8^k}\eq\l(\f {8(8-4)}{p^a}\r)=\l(\f 2{p^a}\r)\ (\mo\ p^2).$$

Clearly $q=p-(\f{4-16}p)=p-(\f p3)$ is divisible by 3 and the two roots of the equation $x^2-4x+16=0$ are
$$2+2\sqrt{-3}=-4\omega^2\ \t{and}\ 2-2\sqrt{-3}=-4\omega,$$
where $\omega=(-1+\sqrt{-3})/2$ is a primitive cubic root of unity.
Thus
$$u_q(4,16)=\f{(-4\omega^2)^q-(-4\omega)^q}{4\sqrt{-3}}=0$$
since $3\mid q$. Applying (1.5) with $m=16$ we get
$$\sum_{k=0}^{(p^a-1)/2}\f{\bi{2k}k}{16^k}\eq\l(\f {16(16-4)}{p^a}\r)=\l(\f 3{p^a}\r)\ (\mo\ p^2).$$

The proof of Corollary 1.1 is now complete. \qed
\medskip

\noindent{\it Proof of Corollary 1.2}. Set $n=(p-1)/2$. Then
$$\align&\sum_{k=1}^n\f{\bi{2k}k}{16^k}\l(\f1{2k-1}+\f1{(2k-1)^2}\r)
\\=&\sum_{k=1}^n\f{2\bi{2k-1}k}{16^k}\cdot\f{2k}{(2k-1)^2}
=\f14\sum_{k=0}^{n-1}\f{\bi{2k}k}{(2k+1)16^k}
\eq0\pmod{p^2}
\endalign$$
with the help of [Su11b, (1.4)].

 Observe that
$$\align \sum_{k=1}^n\f{\bi{2k}k}{(2k-1)16^k}=&\sum_{k=1}^n\f{2\bi{2k-1}k}{(2k-1)16^k}=2\sum_{k=1}^n\f{\bi{2k-2}{k-1}}{k16^k}
\\=&2\sum_{j=0}^{n-1}\f{C_j}{16^{j+1}}=\f18\sum_{k=0}^n\f{C_k}{16^k}-\f{C_n}{8\times 4^{2n}}.\endalign$$
Also,
$$\f{C_n}{4^{2n}}=\f{\bi{p-1}{(p-1)/2}}{4^{p-1}(p+1)/2}\eq(-1)^{(p-1)/2} 2(1-p)\ \ (\mo\ p^2)$$
in view of Morley's congruence ([Mo])
$$\bi{p-1}{(p-1)/2}\eq(-1)^{(p-1)/2}4^{p-1}\ \ (\mo\ p^3).$$
By (1.6) and (1.8),
$$\sum_{k=0}^n\f{C_k}{16^k}\eq-6\l(\f 3p\r)+8-2p\l(\f{-1}p\r)\ \ (\mo\ p^2).$$
Therefore
$$\align\sum_{k=1}^n\f{\bi{2k}k}{(2k-1)16^k}\eq&\f{8-6(\f 3p)-2p(\f{-1}p)}8-\f{2(1-p)(\f{-1}p)}8
\\=&1-\f{(\f{-1}p)+3(\f 3p)}4=1-\l(\f{-1}p\r)\f{3(\f p3)+1}4\ \ (\mo\ p^2)\endalign$$
and hence
$$\sum_{k=1}^n\f{\bi{2k}k}{(2k-1)^2 16^k}\eq-\sum_{k=1}^n\f{\bi{2k}k}{(2k-1)16^k}\eq
-1+\l(\f{-1}p\r)\f{3(\f p3)+1}4\ \ (\mo\ p^2),$$
which yields (1.9) and its equivalent form (1.10). We are done. \qed

\medskip

\Ack. The author would like to thank the referee for helpful comments.

\enddocument